\documentclass{sigcomm-alternate}
\usepackage{amsfonts,latexsym,amsmath,amstext,amssymb,verbatim,epsfig,psfrag}
\usepackage{colordvi,pstcol}
\usepackage{graphicx,psboxit}
\usepackage{psfrag}

\def\N{\mbox{I\hspace{-.15em}N}}
\def\E{\mbox{I\hspace{-.15em}E}}

\newtheorem{theorem}{Theorem}

\begin{document}
\title{Understanding differential equations through diffusion point of view}

\numberofauthors{1}
\author{
   \alignauthor Dohy Hong\\
   \affaddr{Alcatel-Lucent Bell Labs}\\
   \affaddr{Route de Villejust}\\
   \affaddr{91620 Nozay, France}\\
   \email{\normalsize dohy.hong@alcatel-lucent.com}
}

\date{\today}
\maketitle

\begin{abstract}
In this paper, we propose a new adaptation of the D-iteration algorithm to numerically solve the differential equations. This problem can be reinterpreted in 2D or 3D (or higher dimensions) as a limit of a diffusion process where the boundary or initial conditions are replaced by fluid catalysts. Pre-computing the diffusion process for an elementary catalyst case as a fundamental block of a class of differential equations, we show that the computation efficiency can be greatly improved.
The method can be applied on the class of problems that can be addressed by the Gauss-Seidel iteration, based on the linear approximation of the differential equations.
\end{abstract}
\begin{psfrags}
\section{Introduction}\label{sec:intro}
The iterative methods to solve differential equations based on the linear approximation are very well
studied approaches \cite{Johnson_1987}, \cite{Ascher:1998:CMO:551054}, \cite{Podlubny:395913}, \cite{Gear:1971:NIV:540426},
\cite{Smith_1985}, \cite{Golub1996}, \cite{Saad}.
The approach we propose here (D-iteration) is a new approach initially applied to numerically solve the eigenvector
of the PageRank type equation \cite{dohy}, \cite{d-algo}, \cite{dist-test}, \cite{distributed}, \cite{partition}, 
\cite{revisit}.

The D-iteration, as diffusion based iteration, is an iteration method that can be understood as a column-vector
based iteration as opposed to a row-vector based approach. Jacobi and Gauss-Seidel iterations are good examples of
row-vector based iteration schemes. While our approach can be associated to the {\em diffusion} vision, the
existing ones can be associated to the {\em collection} vision.

In this paper, we are interested in the numerical solution for linear equation:
\begin{eqnarray}\label{eq:le}
A.X &=& B
\end{eqnarray}
where $A$ and $B$ are the matrix and vector associated to the linear approximation of
differential equations with initial conditions or boundary conditions. 

In \cite{diff}, it has been shown how simple adaptations can make the diffusion
approach an interesting candidate as an alternative iterative scheme to numerically solve
differential equations. In this paper, we propose a new approach based on the pre-computation
of the elementary diffusion limit. This limit can be then used for a given class of differential
equations, for instance for 2D and 3D case, or for higher dimension.

In Section \ref{sec:heat}, we introduce the 2D problem formulation.
In Section \ref{sec:catal}, we define the notion of catalyst position
and elementary solution. Section \ref{sec:algo} describes the algorithm
with the use of the elementary solution. Finally, Section \ref{sec:eval} gives
an illustration of the application and an evaluation of the run time gain.

\section{From the heat equation}\label{sec:heat}

A typical linearized equation of the stationary heat equation in 2D is of the form:
\begin{align*}
T(n,m) &= \frac{1}{4} \left( T(n-1,m)+T(n,m-1)\right.\\
&\qquad \left.+T(n+1,m)+T(n,m+1)\right)
\end{align*}
which can be obtained by the discretization of the Laplacian operator in Cartesian coordinates:
\begin{eqnarray}
\Delta T(x,y) &=& \frac{\partial^2 T}{\partial x^2} + \frac{\partial^2 T}{\partial y^2} = 0.
\end{eqnarray}
inside the surface $\Omega$ (for instance, $\Omega = [0,L_x]\times[0,L_y]$).
Then additive terms appear for the initial or boundary conditions (Dirichlet) on the
frontier $\partial\Omega$ (for instance, for $x=0$ or $y=0$ etc).

More generally, we consider here linear equations of the form:
\begin{align}
T(n,m) &= \frac{1}{4} \left( T(n-1,m)+T(n,m-1)\right.\\
&\qquad \left.+T(n+1,m)+T(n,m+1)\right)\label{eq:lin}
\end{align}
when $(n,m) \in \Omega^0 = \Omega \setminus \partial\Omega$ and
\begin{align}
T(n,m) &= g(n,m).\label{eq:bc}
\end{align}
when $(n,m) \in \partial\Omega$. For the general case, $\partial\Omega$ is not
necessarily only the boundary of $\Omega$ (at least for discrete formulation, its asymptotic
limit to the initial continuous problem is another story)
and we may add any set of points included in $\Omega$.

We recall that the D-iteration requires updating two vectors (cf. \cite{d-algo}): the fluid vector $F$ and the history
vector $H$ instead of a single vector for the Gauss-Seidel. The above equation can be solved or
by iterating the equation \eqref{eq:lin} on $\Omega^0$ or by applying the diffusion process
associated to the D-iteration.

\section{Elementary solution}\label{sec:catal}

\subsection{Diffusion on 1D}
For the purpose of illustration, let's consider the 1D case with
the following equations (from differential equation of order two):

\begin{align}
T(n) &= \frac{1}{2} \left( T(n-1)+T(n+1)\right)\label{eq:lin1d}
\end{align}
when $n \in \Omega^0$ and
\begin{align}
T(n) &= g(n).\label{eq:bc1d}
\end{align}
when $n \in \partial\Omega$.

The solution can be found by iteration of equation \eqref{eq:lin1d}
on $\Omega^0$ with boundary condition:
\begin{align*}
T(n) &= g(n) \mbox{ when } n \in \partial\Omega\\
T(n) &= 0 \mbox{ when } n \notin \partial\Omega.
\end{align*}

This is exactly the Gauss-Seidel iteration if we always use the latest
updated values of $T$ and the limit is the piecewise linear function
joining the points $(n,g(n))$ at $n\in \partial\Omega$.

The diffusion based iteration would apply the equations:
\begin{align}
H(n) &+= F(n); \mbox{ } F(n) = 0,\\
F(n+1) &+= \frac{1}{2} F(n), \mbox{ }
F(n-1) += \frac{1}{2} F(n)\label{eq:lin1d-diff}.
\end{align}

\subsection{Catalyst position}
We first introduce the notion of catalyst node:
a node (position $n$) is a catalyst if it diffuses once its initial fluid value
$F(n)$, then behaves as a black hole, i.e., it absorbs all fluids it receives without
retransmitting them to its neighbour positions (so that we have $H(n)$ is constant to $F(n)$
after its diffusion).
 
\begin{theorem}
If we associate $\partial\Omega$ to catalyst positions with
the boundary conditions \eqref{eq:bc1d} defining their initial fluid values, the
D-iteration diffusion's limit is equal to the limit of the iteration of the equation \ref{eq:lin1d}.
\end{theorem}
\proof
This theorem is a direct consequence of the equation on $H$ which does exactly
\eqref{eq:lin1d} (cf. \cite{revisit}).

To solve efficiently such an iteration scheme (of course, the interest is
for 2D/3D or more complex mix of differential equations of higher order in 2D/3D),
we introduce the notion of elementary solution for the fluid diffusion process:
the elementary solution is the limit of the diffusion process when we
put a fluid 1 at position $(0,0)$ and with boundary condition zero
on the frontier.
If we impose that $(0,0)$ is a catalyst, the solution can be obtained from
the elementary solution by re-normalizing all values by $1/(1-H(0,1))$. 

We call elementary catalyst when we have $\partial\Omega = \{0,N\}$ and
$g(0) = 1, g(N)=0$ with $\Omega = [0,N]$.

We can compare this idea to the use of the Green's function to solve linear differential
equations. The difference is that the function we define here is much easier to compute
and not dependent on the boundary shape.

Then let's consider the limit of \eqref{eq:lin1d} when we have an elementary
catalyst: by symmetry, we can just explore the space $\N$:
\begin{itemize}
\item at the first iteration, the half of 1 is sent to the position 1;
\item at the second iteration, half of $1/2$ ($1/4$) is sent to positions 2 and 0:
  0 is a catalyst, therefore the $1/4$ it receives disappears;
\item at the thirst iteration, $1/8$ is sent to 1 and 3, and so on.
\end{itemize}

\begin{figure}[htbp]
\centering
\includegraphics[width=6cm]{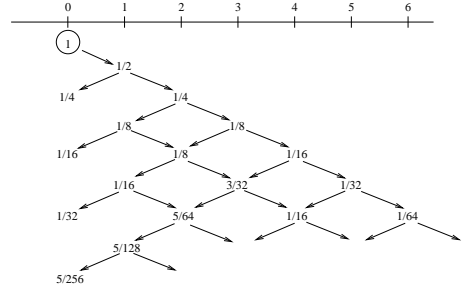}
\caption{Elementary catalyst on 1D: diffusion on $\N$.}
\label{fig:cata1d-fig}
\end{figure}

If we iterates the diffusion, at the limit we find a discrete function which
can be used to solve all equations of type \eqref{eq:lin1d} (see Section \ref{sec:algo}). 
What's interesting is that when we don't put the boundary $N$,
we can prove that the amount of fluid that disappeared in 
the {\em black hole} is associated to the limit of hypergeometric series defined by:
$$
x_{n+1} = \frac{2n+1}{2n+4}x_n
$$
and converges to $1/2$ (use of Gamma function, \cite{de1937fonction}).
Since $1/2$ is the fluid sent to the direction $\N$, that means that all fluid are
finally absorbed by the {\em black hole} and that each value $H(n)$ is convergent
(non-decreasing and bounded). What's even more interesting is that at the limit
we have $H(n)=1$ for all $n\in \N$ (proof by induction from what's received at $0$),
which means that all positions will receive and send exactly 1 fluid: this is a consequence
of the fact that the associated random walk is recurrent (only for 1D and 2D).
However the convergence speed to the limit of this process (by diffusion or by
collection) is very slow (for large $n$). This is why it may be interesting to pre-compute
those iterations once, but not up to its limit (which is the constant function), but
with a boundary condition at $N$.

In the equation we consider, we'll always have the frontier of $\Omega$ that
are all catalysts, therefore all fluid reaching the border of $\Omega$ will
disappear and this guarantees a faster convergence compared to the case without
boundary. Since the limit is the limit of the diffusion iterations when all
catalysts have injected exactly $g(n)$ and when $|F(n)| = 0$, we can use
the catalyst limit we pre-compute on a finite set 
as a common block for faster diffusion (see Section \ref{sec:algo}).

\begin{figure}[htbp]
\centering
\includegraphics[angle=-90,width=6cm]{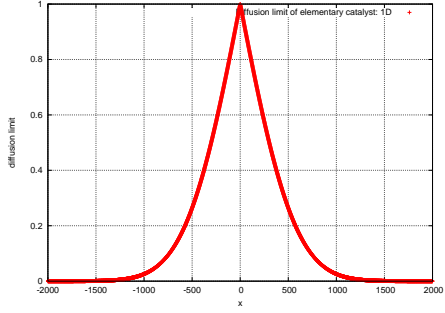}
\caption{Diffusion of the 1D elementary catalyst on $[-2000, 2000]$.}
\label{fig:cata1d}
\end{figure}

Figure \ref{fig:cata1d} shows the limit
function we obtained after $10^5\times 2000$ iterations ($10^5$ cycles over $[1,..2000]$):
note that this is far from linear function! As mentioned above, 
the limit is in this case the piecewise linear function that joins $(-2000,0)$ to $(0,1)$
and $(0,1)$ to $(2000,0)$.

From the linear algebra point of view, this means that we have a matrix associated to the equation
\eqref{eq:lin1d} with initial condition $g(0)=1$ and that its spectral radius is strictly less
than 1 thanks to the boundary condition.

We can also interpret this limit (after normalization) as the average sojourn time of a random walker
when starting from position $0$ and to which we forbid to return to $0$ or reach the boundary position
$N$.

Remark that, if $0$ was not a catalyst and sends back the fluid it receives, we end up with
unbounded quantities (for $H$) when the boundary $N$ is not set,
since the initial fluid 1 never disappears (null-recurrent random walk, stating
from $0$). 

Finally, note that in the general case, we will have equation of the form
(instead of Equation \eqref{eq:lin1d}):
\begin{align}\label{eq:lin1d-gen}
T(n) &= \alpha T(n-1)+\beta T(n+1).
\end{align}
An example is given in the next section.
We may have also equations where $\alpha$ and $\beta$ may depends
on $n$. An example of this is the discretization of the Laplacian in polar coordinates,
invariant by rotation:
\begin{align}\label{eq:lin1d-lap-pol}
\Delta T &=  \frac{\partial^2 T}{\partial r^2}+\frac{1}{r}\frac{\partial T}{\partial r},
\end{align}
which gives:
\begin{align}\label{eq:lin1d-lap-pol-dis}
T(n) &=  \frac{2n+1}{4n}T(n+1)+\frac{2n-1}{4n}T(n-1).
\end{align}

\subsection{Application on 1D}
Consider a differential equation of 2nd order:
\begin{align}
& y{''}(x) + \alpha y{'}(x) + \beta y(x) = f(x)
\end{align}
with boundary condition $y(0)$ and $y(L_x)$.
The naive corresponding iteration scheme (from discretization) is:
\begin{align}
y_n =& \frac{1}{2+\alpha \epsilon -\beta \epsilon^2} \times\\
 &\left( (1+\alpha \epsilon)y_{n+1} + y_{n-1} - \epsilon^2 f_n\right)
\end{align}
or
\begin{align}\label{eq:1D2O}
y_n =& \frac{1}{(2-\beta \epsilon^2)} \times\\
 &\left( (1+\alpha \epsilon/2)y_{n+1} + (1-\alpha \epsilon/2)y_{n-1} - \epsilon^2 f_n\right)
\end{align}
for increment $\epsilon$ ($y_n=y(n\epsilon)$ and $f_n = f(n\epsilon)$).
To solve this iteration scheme, we can:
\begin{itemize}
\item solve the elementary catalyst pre-diffusion;
\item then use the catalyst pre-diffusion to:
  \begin{itemize}
  \item for each position $n$, diffuse $\epsilon^2 f_n$;
  \item solve the diffusion problem with boundary $y(0)$ and $y(L_x)$.
  \end{itemize}
\end{itemize}

For the diffusion of $\epsilon^2 f_n$ we can use a pre-diffusion as for elementary
catalyst but without the constraint of the initial position behaving as a
black hole. Then the elementary catalyst diffusion limit can be obtained by
normalizing all values by $H(0)$.

\subsubsection{Simple example}\label{sec:1Dexm}
Assume we want to solve the equation:
\begin{align*}
& y''(x) = f(x)
\end{align*}
with $f(x) = (-0.99 cos(x) + 0.2 sin(x)) \exp^{-x/10}$ and with boundary
conditions: $y(0) = 1, y(50)=0$ (the solution is $y(x)=cos(x) e^{-x/10}$).

The discretized equation is (from Equation \eqref{eq:1D2O} with $\alpha=\beta=0$):
\begin{align}
y_n =& \frac{1}{2}(y_{n+1} + y_{n-1}) - \frac{\epsilon^2}{2} f_n
\end{align}

Since the catalyst limit of \eqref{eq:lin1d} is the piecewise linear function,
we can solve the equation as follows:
\begin{itemize}
  \item for each position $n$, diffuse $\epsilon^2 g_n/2$, which means adding
    the linear function (multiplied by $L_x = 1/(1-(L_x-1)/L_x)$, because $(L_x-1)/L_x$
    is the quantity that comes back to the diffusion initialization point and here
    there is no {\em black hole} behaviour):
\begin{verbatim}
  for (int i=1; i < Lx-1; ++i){
    transit = - step*step/2 * g(step*i) * (Lx-1);
    for (int j=-i; j < Lx-i; ++j){
      int x = i+j;
      Y[x] += transit * (Lx-1 - abs(j)) / (Lx-1);
    }
  }
\end{verbatim}
(where step is $\epsilon$);
  \item solve the diffusion problem with boundary $y(0)$ and $y(L_t)$:
\begin{verbatim}
  transit = 1.0 - Y[0];
  for (int i=0; i < Lx; ++i){
    Y[i] += transit * (Lx-1 - i) / (Lx-1);
  }
  transit = 0.0 - Y[Lx-1];
  for (int i=0; i < Lx; ++i){
    Y[i] += transit * i / (Lx-1);
  }
\end{verbatim}
\end{itemize}

For any $\epsilon$ we obtain directly the limit
as a superposition of fluids diffusion from $f_n$ and the boundary conditions and
there is no need to do any iterations thanks to the explicit form of the catalyst limit.
The results are shown on Figure \ref{fig:1D2O}.

\begin{figure}[htbp]
\centering
\includegraphics[angle=-90,width=\linewidth]{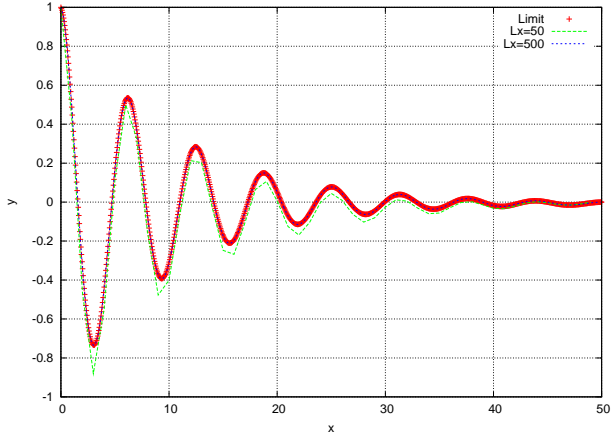}
\caption{Example of 1D.}
\label{fig:1D2O}
\end{figure}



The explicit theoretical formulation of the algorithm in this example is:
\begin{align}\label{eq:1D-var-finite}
y(\alpha L_x) &= (1-\alpha)y(0)+\alpha y(L_x)\\
&+\frac{L_x^2}{2N^2}\sum_{i=1}^{N-1}\left(2\alpha i -i+|N/L_x\alpha-i|-\alpha N\right) f(L_x/N i)
\end{align}
which in the limit ($N\to\infty$) is for $L_x=1$:
\begin{align}\label{eq:1D-var-lim}
y(x) &= (1-x)y(0) + x y(1)\\
&+\frac{1}{2}\int_0^1\left(2x t -t-x+|x-t|\right) y''(t) dt.
\end{align}

This formulation can be directly solved if we look for an expression of
$y(x)$ as a function of $x, y(0), y(1)$ and of the integral of the form $\int u(t,x)y''(t) dt$.
We'll see how such a formulation can be generalized. In this case, the diffusion approach
is equivalent to the usual discretization of the integral in the equation \eqref{eq:1D-var-lim}.

\subsection{Diffusion on 2D}
In this section, we consider linear equations associated to:
\begin{align}
T(n,m) &= \frac{1}{4} \left( T(n-1,m)+T(n,m-1)\right.\\
&\qquad \left.+T(n+1,m)+T(n,m+1)\right)\label{eq:lin2d}
\end{align}
when $(n,m) \in \Omega^0$ and
\begin{align}
T(n,m) &= g(n,m).\label{eq:bc2d}
\end{align}
when $(n,m) \in \partial\Omega$.

For 2D, we consider as for 1D, the diffusion limit of the elementary catalyst at position $(0,0)$.
Figure \ref{fig:cata2d} shows the limit function we obtained.

\begin{figure}[htbp]
\centering
\includegraphics[angle=-90,width=\linewidth]{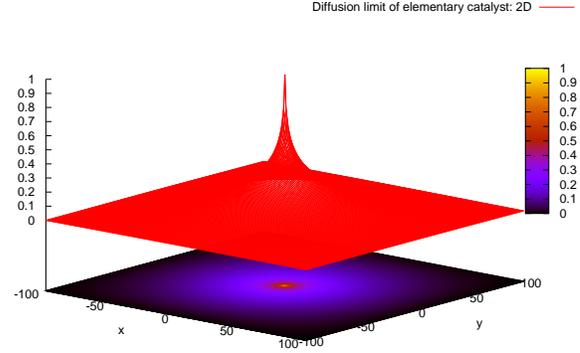}
\caption{Limit of the 2D elementary catalyst on $[-100, 100]\times [-100,100]$.}
\label{fig:cata2d}
\end{figure}

The computation of this limit on a large space is computation costly.
Using the rotation invariant polar Laplacian equation, we can in fact find
the explicit solution, which is of the form:
$$
C - \frac{\log(r)}{B}.
$$
This function has a singularity at 0 (because the limit to the continuous case must be a density or a measure).
An empirical interesting candidate to approximate the diffusion limit of the elementary catalyst
in 2D is:
$$
T(n) = \alpha \left(1 - \frac{\log(n)}{\log(L_r)}\right),
$$
where $\alpha$ has been evaluate from the explicit diffusion iterations.
$T(n)$ is then such that $T(0)=1$ by definition, $T(1)=\alpha$ and $T(L_r)=0$.
In fact, we found that it was better to use the iteration of Equation \eqref{eq:lin1d-lap-pol-dis}
which was close to the above close formula, except the tails.
However, those limits are in polar coordinates, which introduce a bias
when applied on the Cartesian coordinates and which we can not eliminate
(because diffusion on grid is not rotation invariant!).  
In this paper (Section \ref{sec:eval}), we used first the iteration of Equation \eqref{eq:lin1d-lap-pol-dis}
(which is 1D, so very fast), then from this we defined the starting point of the iteration
on 2D, using also the symmetry of 2D (computation on 1/8 of the plane), which accelerated
the full naive 2D scheme iteration by factor 10-50.

Another alternative is to apply the ideas of Section \ref{sec:algo} during the
pre-computation: after iterations on a smaller space (for instance, $N'=N/2$), 
we can save the results $H^0, F^0$, 
then we can replace the elementary diffusion
by directly copying the results of $N$ iterations as a block.
This is interesting to gain an order of precision quickly, exploiting the fact
that $F^0$ has fluids concentrated at the border. However, after one block copy
operations, we find again a configuration where the fluids are spread more
uniformly. Optimizing the pre-computation phase is an independent problem which we don't
analyse further here.

\section{Algorithm}\label{sec:algo}
We consider the 2D problem 
$\Delta y = f$ on $\Omega$ with boundary condition $g$ on $\partial\Omega$ for
illustration. The method should be easily extended to a much general
linear operator associated to other differential equations.

We assume the elementary catalyst's limit is pre-computed on a finite set
$[-L_x,L_x]\times[-L_y,L_y]$ with boundary condition $g(0,0) = 1$ and
$g(x,y)=0$ if $|x|= L_x$ or $|y| = L_y$. For the sake of simplicity, 
we will consider $\Omega$ of the form $[0,L_x]\times[0,L_y]$ (if not,
we can choose $L_x$ the maximal x-distance between two points of $\Omega$
and similarly for $L_y$).
For the practical computation, we iterate the D-iteration until the
remaining fluid $|F|$ is below the targeted error. Then, we store
in a file the last states on $H$ and $F$ (in the following denoted
$H^0$ and $F^0$).

Then we apply the following process:
\begin{itemize}
\item load the above results $H^0$ and $F^0$;
\item define a new variable $H[L_x][L_y]$ and $F[L_x][L_y]$ (initialized to 0);
\item set the initial fluid $F$ equal to $f$ in $\Omega^0$;
\item diffuse $F$ on $\Omega$ (including $\partial\Omega$): 
  here, diffusion means adding $H^0$ and $F^0$ on $H$ and $F$ respectively at translated
  position (by $x$ and $y$);
\item diffuse fluid $g(x,y)-H[x][y]$ on $\partial\Omega$ (choose positions where
  $|g(x,y)-H[x][y]|$ is the largest or above a certain threshold;
  here, diffusion means adding $H^0$ and $F^0$ on $H$ and $F$ respectively at translated
  position (by $x$ and $y$):
\begin{verbatim}
 Diffusion of "g(x,y)-H[x][y]" :
  for (int x=0; x < Lx; ++x){
    for (int y=0; y < Ly; ++y){
      if ( bound[x][y] ){// boundary position
        transit = g[x][y] - H[x][y];
        if ( abs(transit) > Thresh_ ){
          for (int i=0; i < n_x; i++){
            for (int j=0; j < n_y; j++){
              H[i][j] += 
           transit*H0[abs(i-x)][abs(j-y)];
              F[i][j] += 
           transit*F0[abs(i-x)][abs(j-y)];
            }
          }
        }
      }
    }
  }
\end{verbatim}
\item the previous step is repeated until the threshold is below the targeted error;
\item if required, we may also diffuse fluid $F$ which are above a given threshold
  (we may also decide not to use $F$ at all), because
  as far as we keep $F$ and $H$, all operations are invertible in the sense that we 
  can inject the surplus or the deficit fluid to make the exact convergence in any order.
\end{itemize}

The numerical solution to our problem is then given by $H$.
If $H$ is exactly equal to the boundary condition $g$ on $\partial\Omega$, then
$H$ on $\Omega$ is the exact limit.

As for the 1D case, we can express this approach by the projection method
where the elementary catalyst limit serves as a unique base.
It can be also understood as an application of calculus of variations or
a Lagrangian approach. Let's call $\phi$ the limit of the elementary catalyst
(for instance on a square surface that's bigger than $\Omega$).
Then, we can rewrite the algorithm under the form:
\begin{align}\label{eq:2D-var-finite}
y(x) &\sim y(n,m)\\
y(n,m) &= -\frac{\delta^2}{4}\sum_{i\in\Omega} f(i)\phi_i(x)\\
&+\sum_{x_b\in\partial\Omega}\left(y(x_b) - \alpha(x_b) - \frac{\delta^2}{4}\sum_{i\in\Omega}g(i)\phi_i(x_b)\right)\tilde{\phi}_{x_b}(x),
\end{align}
where $x=(\delta n, \delta m)\in\Omega$ (regular grid of $\delta$), $\phi_i(x)$ the value of
$\phi$ at point $x$ when the origin is set at $i$ and $\tilde{\phi} = \phi/(1-\phi(0,1))$,
and $\alpha(x_b)$ is term expressing all diffusion received from other boundary condition
related diffusion.
Our approach can be understood as an iterative approach to find the coefficients
$\alpha(x_b)$.

Its limit (if existence) for $\delta\to 0$ can be formulated as:
\begin{align}\label{eq:2D-var-lim}
y(x) &= \int_{\partial\Omega}(y(x_b)-\alpha(x_b))\tilde{\phi}(x_b-x)dx_b\\
&-\frac{1}{4}\int_{\Omega} \Delta y(t)\phi(x-t)dt\\
&-\frac{1}{4}\int_{\partial\Omega}\int_{\Omega} \Delta y(t)\phi(x_b-t)\tilde{\phi}(x_b-x)d t d x_b
\end{align}
where the second term comes from the diffusion of fluid inside $\Omega$
and the two other from the correction for the boundary conditions.
This formula assume in particular that we have a limit of $\phi$ when $N$ goes to infinity.
We can interpret $\phi(n,m)$ as the probability for 2D random walk to reach
$(n,m)$ before touching the boundary starting from $(0,0)$.
When $N$ goes to infinity, the random walk tends to the 2D Brownian motion
and $\phi(x)$ in the continuous space is the probability that
from $(0,0)$ we reach $[x,x+dx]\times[y,y+dy]$ before the boundary is
touched.

In a particular case when $f=0$, they is a very nice theory of probability
which shows that $y$ is given by an explicit integration formula
(\cite{stroock2005introduction, karatzas1991brownian, dellacherie1975probabilites}):
$$
y(x) = h(x) = \E_x[g(B(T))]
$$
where $B$ is the Brownian motion, $T$ the stopping time when the boundary is touched.
If the boundary is a sphere, we have a more explicit formula of the form:
$$
y(x) = h(x) = \int_{S^{d-1}} \frac{1-|x|^2}{|x-y|^d} g(y) \sigma_d(dy).
$$

From the diffusion point of view, we can understand why with the sphere
we can have a simpler formula: the diffusion from one point of the sphere to all
others points of the sphere follows exactly the same process, meaning that
in our approach the terms $\alpha(x_b)$ can be eliminated if
$\phi$ is associated to this diffusion model.

Our approach can be understood as an explicit practical solution, not
only in presence of $f$, but also for a general operators (so not only harmonic
functions) associated to the differential equations, using a specific choice
of $\phi$.
When the diffusion operator is not symmetrical in the four directions,
the very nice theory of harmonic function does no more apply.
However, the idea of exploiting the pre-diffusion ($\phi$) can be also compared to
the use of the Green's function $G(x,s)$ (when it is known!) and express the solution as:
$$
y(x) = \int G(x,s) f(s) ds 
$$
 
But while this is an exact solution, the computation of the Green's function
may be even more complex than solving directly by an iterative scheme in
a general case.

Note that our algorithm has no guarantee of convergence (on $\alpha(x_b)$).
We hope address this point in a future paper, if such a consideration is not
already proposed in the past.

\subsection{Error estimate}
The distance to the limit can be estimated from
$$
r = \sum_{x,y\in \Omega^0} |F[x][y]| + \sum_{x,y\in\partial\Omega} |g(x,y)-H[x][y]|.
$$
The first component of $r$ is the residual fluid resulting from the diffusion
by catalysts and the second component is the surplus or the deficit fluid
that are injected to $\Omega$.
If $r=0$, $H$ is the exact limit of the problem.

\section{Evaluation}\label{sec:eval}
\subsection{Convergence comparison}
For the evaluation purpose, we considered the following (too simple!) scenario:
\begin{itemize}
\item S1: a 2D diffusion problem ($f=0$) with $L_x=L_y$ and
  boundary condition on the border: $g(x,y) = 100$. The solution of
  this problem is obviously a constant function equal to $100$ on
  every point of $\Omega$.
\end{itemize}

The results are shown on Table \ref{tab:eval0}: we used 2 Linux laptop:
Intel(R) Core(TM)2 CPU, U7600, 1.20GHz, cache size 2048 KB (Linux1, $g++-4.4$)
and
Intel(R) Core(TM) i5 CPU, M560, 2.67GHz, cache size 3072 KB (Linux2, $g++-4.6$).
The pre-computation of the elementary catalyst on $[-L_x, L_x]\times [-L_y, L_y]$
has been done for a given target error (target, on the remaining fluid); the
runtime for this pre-computation is given by pre-comp. The results have been
saved in a simple ASCII file, its loading time is given by Init.
We observed that the limitation of the error of our approach was
about $10^{-5}$ (which means for $g(x)=100$ a relative precision of
$10^{-7}$), resulting probably from the double precision (about $10^{-15})$) we have on
$F^0$ (relatively to $H^0$).
Through, this school case, we just want to illustrate the potential of our
approach.

\begin{table}
\begin{center}
\begin{tabular}{|l|cccc|}
\hline
         & GS         &           & DI &  \\
$L_x$    & 100        & 200       & 100         & 200 \\
\hline
Linux 1  & & & &\\
Pre-comp & x          & x         & 1.2          & 20\\
target   & x          & x         & $1e^{-3}$  & $1e^{-3}$\\
Init     & x          & x         & 0.1         & 0.3 \\
\hline
\hline
error    & 1.0        & 1.0       & 0.8         & 1.0  \\
error2   & $1e^{-3}$  & $2e^{-4}$ & 200         & 400\\
time     & 0.6        & 10        & 0.02        & 0.12 \\
gain     & 1          & 1         & $\times$30  & $\times$80 \\
\hline
error    & 0.1        & 0.1       & 0.1         & 0.09 \\
error2   & $1e^{-4}$  & $2e^{-5}$ & 30          & 53\\
time     & 0.9        & 15        & 0.07        & 0.5 \\
gain     & 1          & 1         & $\times$13  & $\times$30 \\
\hline
\hline
\end{tabular}\caption{Comparison of computation cost. Pre-comp: pre-computation time of the elementary catalyst. target: target error on the remaining fluid for pre-computation. Init: initialization time. error: distance to the limit. error2: maximum increment of the last iteration for GS, $r$ for DI. $L_x=100, 200$.}\label{tab:eval0}
\end{center}
\end{table}

\begin{table}
\begin{center}
\begin{tabular}{|l|cccc|}
\hline
         & GS         &           & DI &  \\
$L_x$    & 300        & 400       & 300         & 400 \\
\hline
Linux 1  & & & &\\
Pre-comp & x          & x         & 100         & 330\\
target   & x          & x         & $1e^{-3}$   & $1e^{-3}$\\
Init     & x          & x         & 0.7         & 1.2 \\
\hline
\hline
error    & 1.0        & 1.0       & 1.0         & 1.0 \\
error2   & $10^{-4}$  & $5e^{-5}$ & 560         & 700\\
time     & 53         & 170       & 0.5         & 1.3 \\
gain     & 1          & 1         & $\times$100 & $\times$130 \\
\hline
error    & 0.1        & 0.1       & 0.1         & 0.1  \\
error2   & $1e^{-5}$  & $5e^{-6}$ & 77          & 100\\
time     & 80         & 250       & 2.2          & 6\\
gain     & 1          & 1         & $\times$35   & $\times$42 \\
\hline
\end{tabular}\caption{Comparison of computation cost: $L_x=300, 400$.}\label{tab:eval1}
\end{center}
\end{table}

\begin{table}
\begin{center}
\begin{tabular}{|l|cccc|}
\hline
         & GS         &        & DI          &\\
$L_x$    & 1000       & 2000   & 1000        & 2000\\
\hline
Linux 1  & & & &\\
Pre-comp & x          & x      & x           & x\\
target   & x          & x      & x           & x\\
Init     & x          & x      & 8           & 15\\
\hline
error    & 1.0        & 1.0    & 1.4         & 0.7\\
error2   & $1e^{-5}$  & $2e^{-6}$& 2000      & 2000\\
time     & 6500       & 105000 & 10          & 90 \\
gain     & 1          & 1      & $\times$650 & $\times$1200 \\
\hline
error    & 0.1        &  0.4   & 0.14        & 0.09\\
error2   & $1e^{-6}$  & $1e^{-6}$& 260         & 270\\
time     & 9500       & 123700 & 53          & 500 \\
gain     & 1          & 1      & $\times$180 & $\times$250 \\
\hline
\hline
Linux 2  & & & &\\
Pre-comp & x          & x      & 6900         & 93000\\
target   & x          & x      & $10^{-3}$    & $10^{-3}$\\
Init     & x          & x      & 3            & 10\\
\hline
error    & 1.0        & 1.0    & 1.0          & 0.7\\
error2   & $1e^{-5}$  & $1e^{-6}$ & 1500      & 2000\\
time     & 2800       & 50000  & 5            & 30\\
gain     & 1          & 1      & $\times$560 & $\times$1600 \\
\hline
error    & 0.1        & 0.1    & 0.11         & 0.09\\
error2   & $1e^{-6}$  & $1e^{-7}$ & 140       & 270\\
time     & 4200       & 55000  & 30           & 150\\
gain     & 1          & 1      & $\times$140 & $\times$370 \\
\hline
\end{tabular}\caption{Comparison of computation cost: $L_x=1000, 2000$. Pre-computation is only done on Linux2.}\label{tab:eval}
\end{center}
\end{table}

\subsection{Stationary heat diffusion in 2D}
Let's consider a simple variant of S1: we set
\begin{itemize}
\item S2: a very simple diffusion problem with $L_x=L_y=2000$ and
  boundary condition on the border: $g(0,y) = 100$ and $g(x,0)=g(x,L_y)=g(L_x,y)=0$. 
\end{itemize}

Results are on Figure \ref{fig:100gs1}, \ref{fig:100gs2}, \ref{fig:100gs4},
\ref{fig:100gs16}, \ref{fig:100di2}, \ref{fig:100di8} and \ref{fig:100di30}: 
for the D-iteration, we use the pre-computation $\phi$
that is generated in the previous section. We can see that with the naive iterative
method, the convergence to the limit may be really slow when a large grid is
considered. The result obtained in 30s with our approach has in this case
a better convergence than with 16 hours with Gauss-Seidel (the gain is reaching
a factor $2000$). But of course, the gain was obtained thanks to the
previous pre-computation $\phi$ which took about 1 day.

\section{Conclusion}\label{sec:conclusion}
In this paper we addressed a first analysis of the potential of the
D-iteration when applied in the context of the numerical solving of
differential equations. We showed that using the regularity of the
diffusion process, we can exploit the idea of the pre-diffusion.
The diffusion approach gives a new way of understanding the differential and
integration associated operator iteration at a fundamental level and
offers a great potential for a very fast numerical computation.
Further exploitation of this will be addressed in a future paper.

\end{psfrags}
\bibliographystyle{abbrv}
\bibliography{sigproc}

\begin{figure}[htbp]
\centering
\includegraphics[angle=-90,width=7cm]{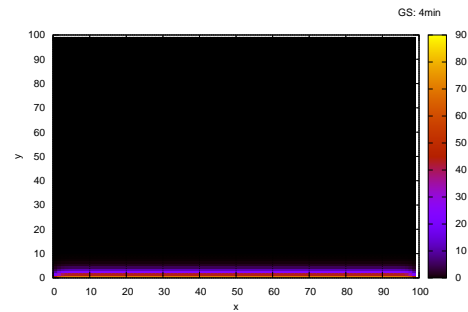}
\caption{Gauss-Seidel. Run time: 4 min.}
\label{fig:100gs4m}
\end{figure}
\begin{figure}[htbp]
\centering
\includegraphics[angle=-90,width=7cm]{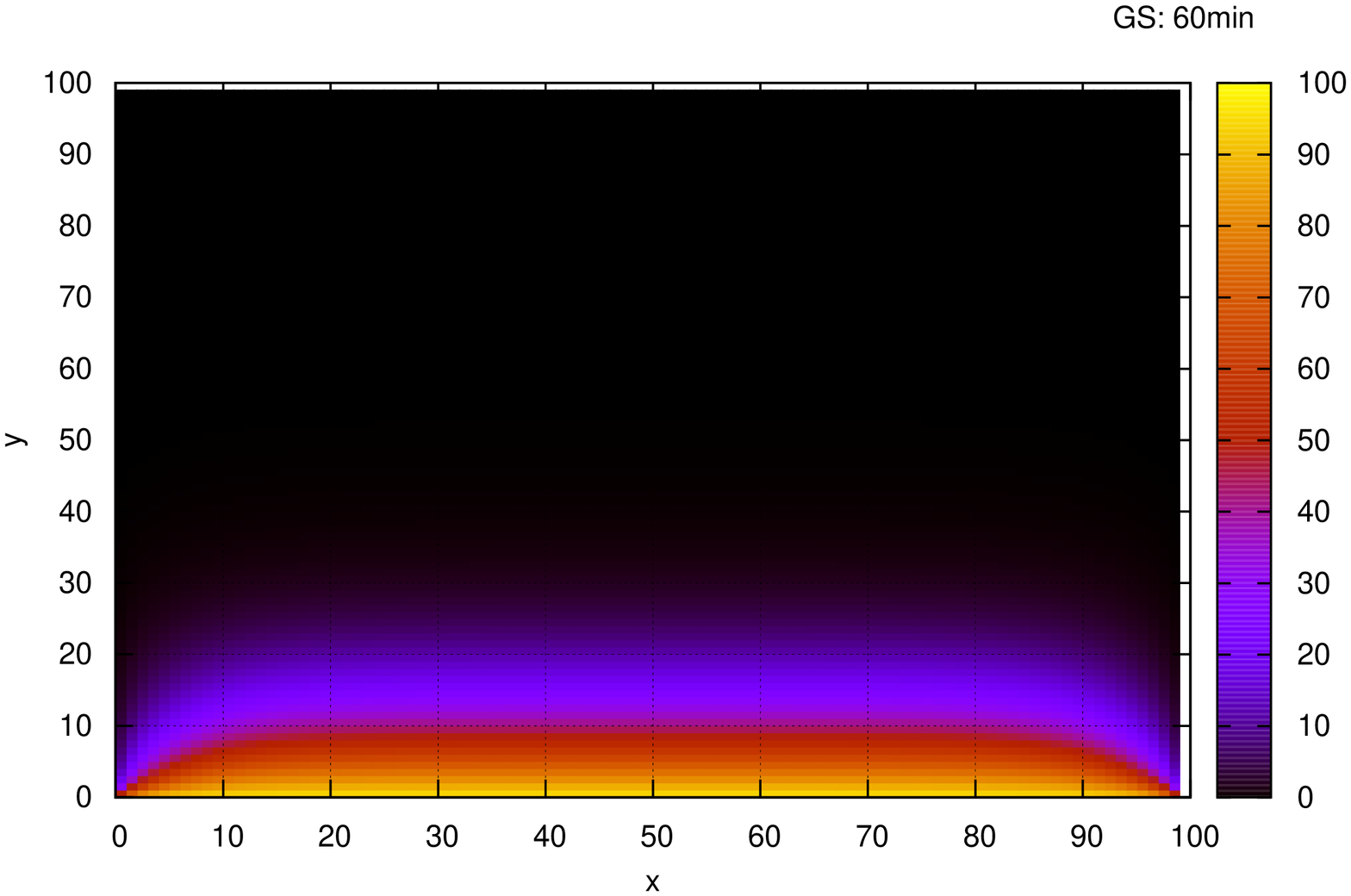}
\caption{Gauss-Seidel. Run time: 1 hour.}
\label{fig:100gs1}
\end{figure}
\begin{figure}[htbp]
\centering
\includegraphics[angle=-90,width=7cm]{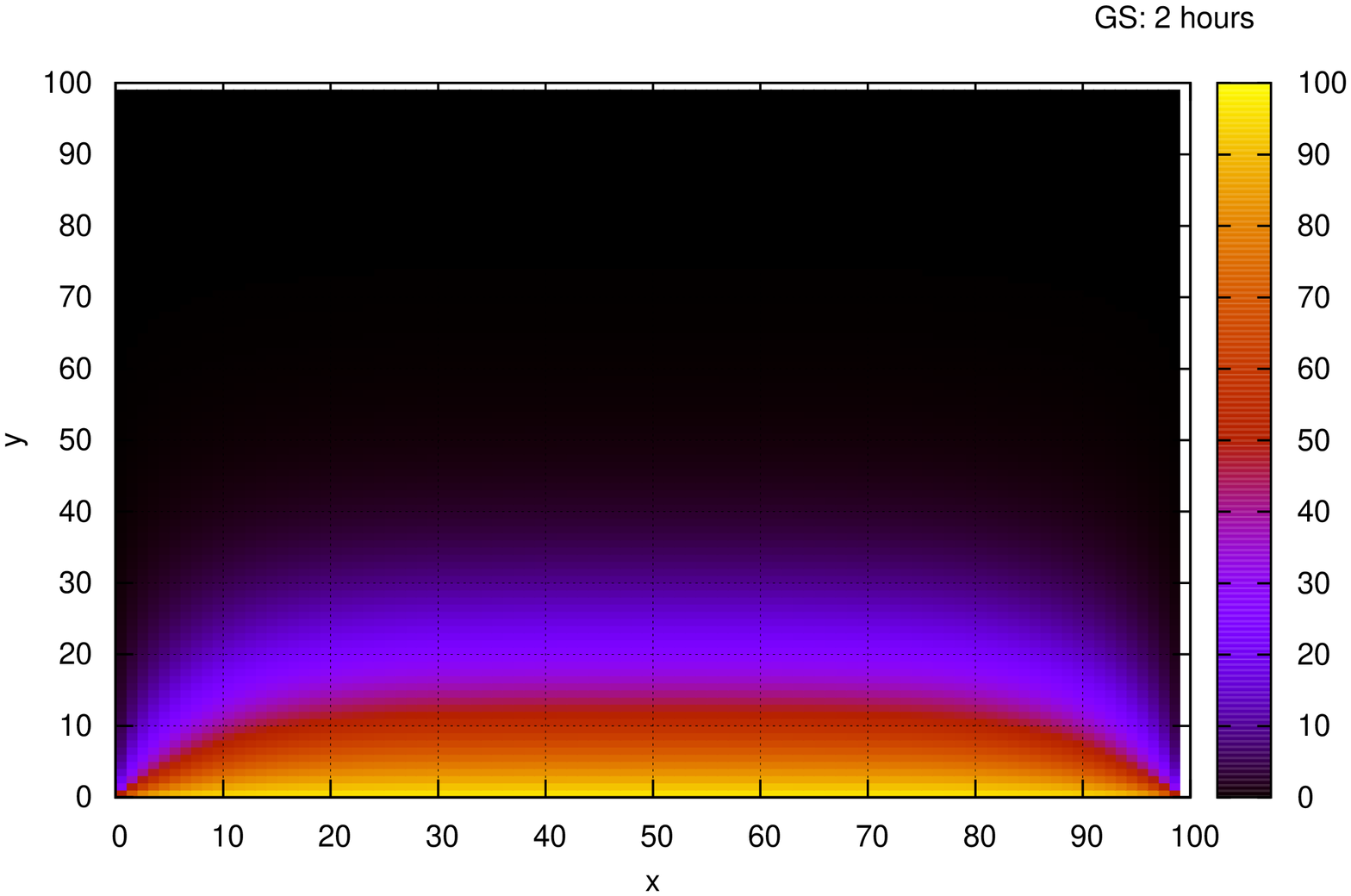}
\caption{Gauss-Seidel. Run time: 2 hours.}
\label{fig:100gs2}
\end{figure}
\begin{figure}[htbp]
\centering
\includegraphics[angle=-90,width=7cm]{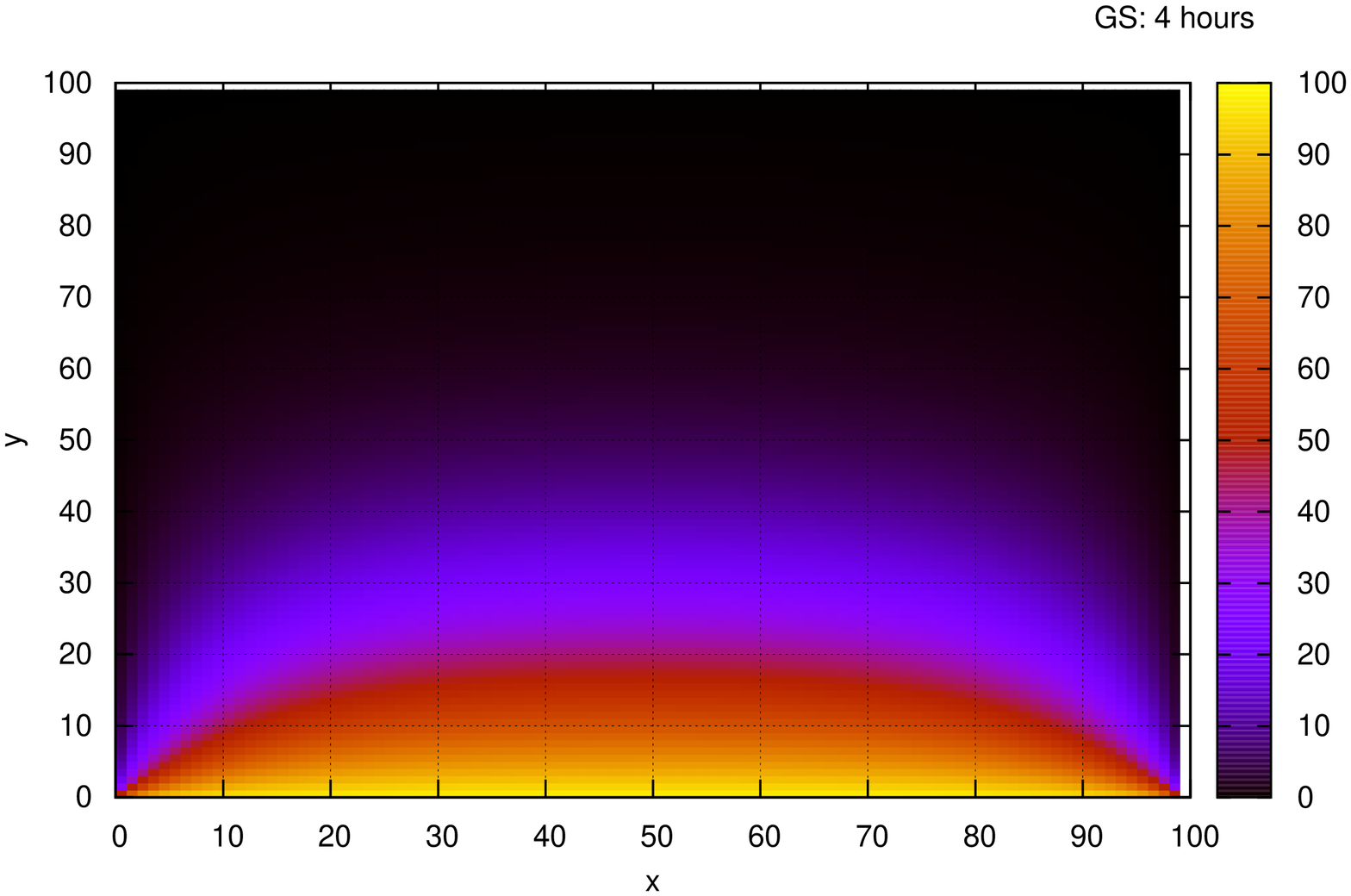}
\caption{Gauss-Seidel. Run time: 4 hours.}
\label{fig:100gs4}
\end{figure}
\begin{figure}[htbp]
\centering
\includegraphics[angle=-90,width=7cm]{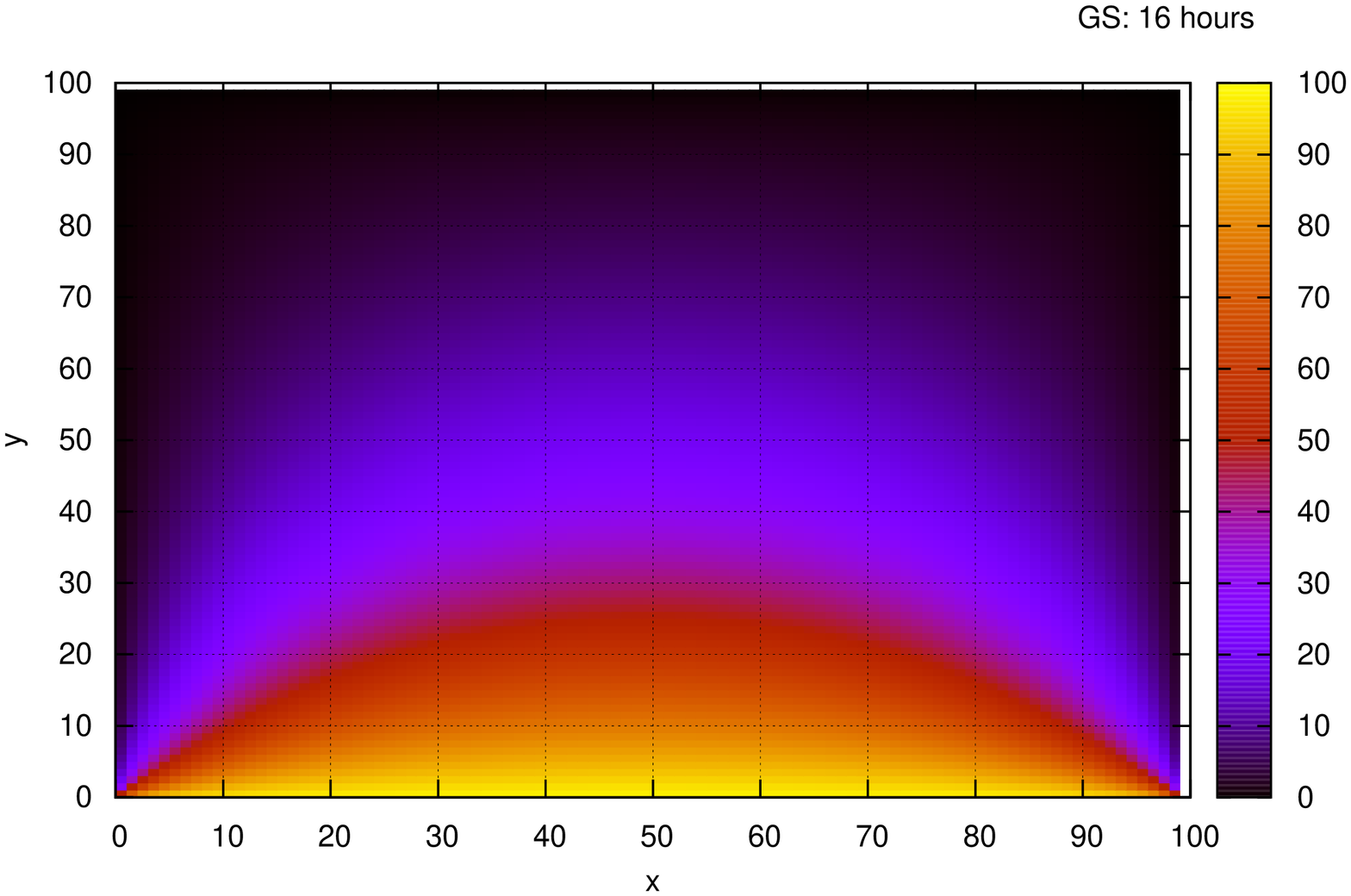}
\caption{Gauss-Seidel. Run time: 16 hours.}
\label{fig:100gs16}
\end{figure}
\begin{figure}[htbp]
\centering
\includegraphics[angle=-90,width=7cm]{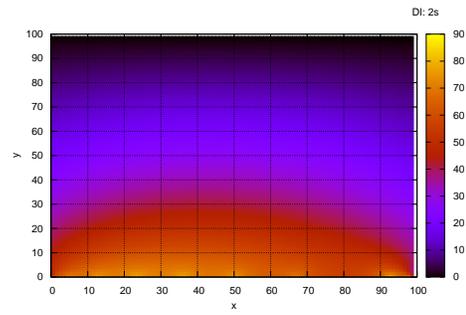}
\caption{D-iteration: 2 s.}
\label{fig:100di2}
\end{figure}
\begin{figure}[htbp]
\centering
\includegraphics[angle=-90,width=7cm]{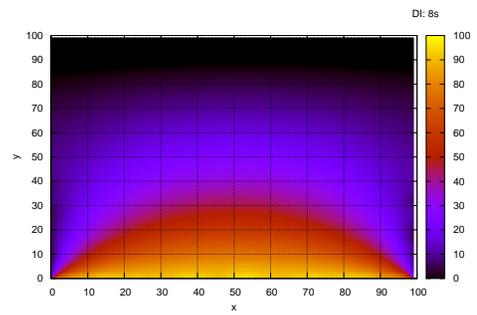}
\caption{D-iteration: 8 s.}
\label{fig:100di8}
\end{figure}
\begin{figure}[htbp]
\centering
\includegraphics[angle=-90,width=7cm]{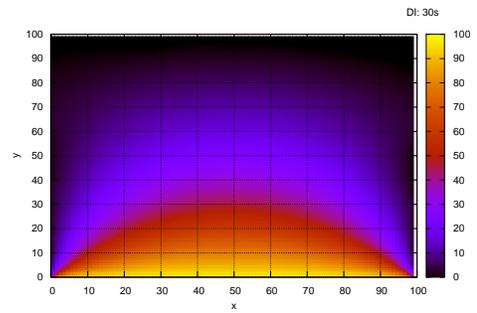}
\caption{D-iteration: 30 s.}
\label{fig:100di30}
\end{figure}

\end{document}